\documentclass[11pt,titlepage]{article}
\usepackage{latexsym}
\usepackage{amsfonts}
\usepackage{amsmath}

\usepackage{amssymb}
\usepackage[margin=1in]{geometry}
\usepackage{graphicx}
\usepackage{tikz}
\usetikzlibrary{arrows, decorations.markings}
\usetikzlibrary{fit}					
\usetikzlibrary{backgrounds}

\newcommand\blackslug{\hbox{\hskip 1pt \vrule width 4pt height 8pt depth 1.5pt
        \hskip 1pt}}
\newcommand\bbox{\hfill \quad \blackslug \medbreak}

\newtheorem{theorem}{}[section]

\newcommand{\Proof}{\noindent{\bf Proof.}\ \ }
\newcommand{\Analysis}{\noindent{\bf Proof of Theorem}\ \ }

\tikzstyle{every node}=[circle, draw, fill=black!50, inner sep=0pt, minimum width=20pt]

\tikzstyle{input}=[circle,
                                    thick,
                                    minimum size=2.2cm,
                                    draw=black!80,
                                    fill=white!20]

\tikzstyle{input2}=[circle,
                                    thick,
                                    minimum size=1.0cm,
                                    draw=black!80,
                                    fill=white!20]

\tikzstyle{matrx}=[rectangle,
                                    thick,
                                    minimum size=1.8cm,
                                    draw=black!80,
                                    fill=white!20]

\tikzstyle{matrx2}=[rectangle,
                                    thick,
                                    minimum size=1.5cm,
                                    draw=black!80,
                                    fill=gray!20]

\tikzstyle{vecArrow} = [thick, decoration={markings,mark=at position
   1 with {\arrow[semithick]{open triangle 60}}},
   double distance=1.8pt, shorten >= 5.5pt,
   preaction = {decorate},
   postaction = {draw,line width=1.4pt, white,shorten >= 4.5pt}]
\tikzstyle{innerWhite} = [semithick, white,line width=1.4pt, shorten >= 4.5pt]

\tikzstyle{background}=[rectangle,
                                                fill=gray!10,
                                                inner sep=0.2cm,
                                                rounded corners=5mm]

\title{On the Erd\H{o}s-Hajnal conjecture for six-vertex tournaments}
\author{
Eli Berger
\thanks{Partially supported by  BSF grant 2006099 and ISF grant 1581/12}\\
University of Haifa\\
Haifa, Israel
\and
Krzysztof Choromanski\\
Google Research\\
New York, NY, USA
\and
Maria Chudnovsky
\thanks{Partially supported by NSF grants DMS-1550991 and BSF grant 2006099}\\
Princeton University\\
Princeton, NJ, USA
}

\date{July 1, 2015; revised \today}

\newcommand{\Keywords}{the Erd\H{o}s-Hajnal conjecture, prime tournaments, galaxies}

\begin{document}
\maketitle
\begin{abstract}
A celebrated unresolved conjecture of Erd\H{o}s and Hajnal states that for every undirected graph $H$ there exists $\epsilon(H)>0$ such 
that every undirected graph on $n$ vertices that does not contain $H$ as an induced subgraph contains a clique or stable set of size at least $n^{\epsilon(H)}$. 
The conjecture has a directed equivalent version stating that for every tournament $H$ there exists $\epsilon(H)>0$ such that every $H$-free
$n$-vertex tournament $T$ contains a transitive subtournament of order at least $n^{\epsilon(H)}$.
We say that a tournament is \textit{prime} if it does not have nontrivial homogeneous sets. So far the conjecture was proved only for
some specific families of prime tournaments (\cite{chorochudber, choromanski2}) and tournaments constructed according to the so-called \textit{substitution procedure}(\cite{alon}).
In particular, recently the conjecture was proved for all five-vertex tournaments (\cite{chorochudber}), but the question about the correctness 
of the conjecture for all six-vertex tournaments remained open. In this paper we prove that all but at most one six-vertex tournament satisfy the Erd\H{o}s-Hajnal conjecture. 
That reduces the six-vertex case to a single tournament.
\end{abstract}

\maketitle {\bf Keywords:} \Keywords

\section{Introduction}

We denote by $|S|$ the size of a set $S$. Let $G$ be a graph. We denote by $V(G)$ the set of its vertices. 
Sometimes instead of writing $|V(G)|$ we will use shorter notation $|G|$. We call $|G|$ the \textit{size of G}. 
We denote by $E(G)$ the set of edges of a graph $G$.
A \textit{clique} in the undirected graph is a set of pairwise adjacent vertices and a \textit{stable set} in the undirected graph is a set of pairwise nonadjacent vertices.
A \textit{tournament} is a directed graph  such that for every pair $v$ and $w$ of vertices, exactly one of the edges $(v,w)$ or $(w,v)$ exists. 
For a tournament $H$ and a vertex $v \in V(H)$ we denote by $H \setminus \{v\}$ the tournament obtained from $H$ by deleting $v$ and all edges incident with it.
We denote by $H^{c}$ the tournament obtained from $H$ by reversing directions of all edges of $H$.
If $(v,w)$ is an edge of the tournament then we say that $v$ is $\textit{adjacent to}$ $w$ (alternatively: $w$ is an \textit{outneighbor} of $v$) 
and $w$ is \textit{adjacent from} $v$ (alternatively: $v$ is an \textit{inneighbor} of $w$). 
For two sets of vertices $V_{1}$, $V_{2}$ of a given tournament $T$ we say that $V_{1}$  is \textit{complete to} $V_{2}$ (or equivalently $V_{2}$ is \textit{complete from} $V_{1}$) if every vertex of $V_{1}$ is adjacent to every vertex of $V_{2}$. We say that a vertex $v$ is complete to/from a set $V$ if $\{v\}$ is complete to/from $V$.
A tournament is \textit{transitive} if it contains no directed cycle. For a set of vertices $V=\{v_{1},v_{2},...,v_{k}\}$ we say that an ordering $(v_{1},v_{2},...,v_{k})$ is \textit{transitive} if $v_{1}$ is adjacent to $v_j$ for every $i < j$.

If a tournament $T$ does not contain some other tournament $H$ as a subtournament then we say that $T$ is $H$-\textit{free}. 

A celebrated unresolved conjecture~ of Erd\H{o}s and Hajnal is as follows:  

\begin{theorem}
\label{EHConun}
For any undirected graph $H$ there exists $\epsilon(H)>0$ such that every $n$-vertex undirected graph that does not contain $H$ as an induced subgraph contains a clique or a stable of size at least $n^{\epsilon(H)}$.
\end{theorem}

In 2001 Alon, Pach and Solymosi proved (\cite{alon}) that Conjecture~\ref{EHConun} has an equivalent directed version, where undirected graphs are replaced by tournaments and cliques and stable sets by transitive subtournaments, as follows:

\begin{theorem}
\label{EHCon}
For any tournament $H$ there exists $\epsilon(H)>0$ such that every $n$-vertex $H$-free $n$-vertex tournament contains a transitive subtournament of size at least $n^{\epsilon(H)}$.
\end{theorem}

If for a graph $H$ there exists $\epsilon(H)>0$ as in \ref{EHCon}, then we say that \textit{$H$ satisfies the Erd\H{o}s-Hajnal conjecture} 
(alternatively: $H$ has the \textit{Erd\H{o}s-Hajnal property}).\\

A set of vertices $S \subseteq V(H)$ of a tournament $H$ is called \textit{homogeneous} if for every $v \in V(H) \backslash S$ the following holds: either for all $w \in S$ we have: $(w,v)$ is an edge or for all $w \in S$ we have: $(v,w)$ is an edge. A homogeneous set $S$ is called \textit{nontrivial} if $|S|>1$ and $S \neq V(H)$. A tournament is called \textit{prime} if it does not have nontrivial homogeneous sets.\\

The following theorem, that is an immediate corollary of the results given in \cite{alon} and applied to tournaments, shows why prime tournaments are important.

\begin{theorem}
\label{nonprime}
If Conjecture~\ref{EHCon} is false then the smallest counterexample is prime.
\end{theorem}

Therefore of interest is studying the Erd\H{o}s-Hajnal property for prime tournaments. 
We need a few more definitions that we borrow from \cite{chorochudber} and put below for the reader's convenience.

For an integer $t$, we call the graph $K_{1,t}$ a {\em star}. Let $S$ be a
star with vertex set $\{c, l_1, \ldots, l_t\}$, where $c$ is adjacent to 
$l_1, \ldots, l_t$. We call $c$ the {\em center of the star}, and
$l_1, \ldots, l_t$ {\em the leaves of the star}. 
Note that in the case $t=1$ we may choose arbitrarily any one of the two vertices to be the center of the star, and the other vertex is then considered to be the leaf. 

Let $\theta=(v_{1},v_{2},...,v_{n})$ be an ordering of the vertex set $V(T)$ of an $n$-vertex tournament $T$.  We say that a vertex $v_{j}$ is \textit{between} two vertices $v_{i},v_{k}$ under $\theta=(v_{1},...,v_{n})$ if $i<j<k$ or $k<j<i$.
An edge $(v_i,v_j)$ is a {\em backward edge under 
$\theta$} if $i>j$.  The \textit{graph of backward edges under $\theta$},
denoted by $B(T, \theta)$, is the undirected graph that has vertex set $V(T)$,
  and   $v_i v_j \in E(B(T, \theta))$ if and only if  $(v_i,v_j)$ or 
$(v_j,v_i)$ is a backward edge of $T$ under $\theta$.

A {\em right star}
in $B(T, \theta)$ is an induced subgraph with vertex set 
$\{v_{i_0}, \ldots, v_{i_t}\}$, such that \\
$B(T,\theta)|\{v_{i_0}, \ldots, v_{i_t}\}$ is a star with center $v_{i_t}$, 
and  $i_t > i_0, \ldots, i_{t-1}$.  In this case we also  
say that $\{v_{i_0}, \ldots, v_{i_t}\}$ is  a right star in $T$.

A {\em left star}
in $B(T, \theta)$ is an induced subgraph with vertex set 
$\{v_{i_0}, \ldots, v_{i_t}\}$, such that \\ 
$B(T,\theta)|\{v_{i_0}, \ldots, v_{i_t}\}$ is a star with center $v_{i_0}$, 
and  $i_0 < i_1, \ldots, i_t$.    In this case we also  
say that $\{v_{i_0}, \ldots, v_{i_t}\}$ is a  left star in $T$.
A {\em star} in $B(T, \theta)$ is a left star or a right star.

Let $H$ be a tournament and assume there exists an ordering $\theta$ of its 
vertices such that every connected component of $B(H, \theta)$ is  either a 
star or a singleton. We call this ordering a \textit{star ordering}. 
If in addition  every star is either a left star or a right star, and no center of a star is between leaves of another star, then the corresponding
ordering is called a \textit{galaxy ordering} and the tournament $H$ is called a \textit{galaxy}.
The main results of \cite{chorochudber} that we will heavily rely on in this paper are:

\begin{theorem}
\label{galaxy-theorem}
Every galaxy has the Erd\H{o}s-Hajnal property.
\end{theorem}

\begin{theorem}
 \label{five-vertex-theorem}
Every tournament $H$ on at most five vertices has the Erd\H{o}s-Hajnal property.
\end{theorem}

We denote by $K_{6}$ the six-vertex tournament with $V(K_{6})=\{v_{1},...,v_{6}\}$ such that under ordering $(v_{1},...,v_{6})$
of its vertices the set of backward edges is: $\{(v_{4},v_{1}),(v_{6},v_{3}), (v_{6},v_{1}), (v_{5},v_{2})\}$. We call this ordering of vertices of $K_{6}$ the \textit{canonical ordering} (Fig.1).

\begin{center}
\begin{tikzpicture}[every path/.style={>=latex},every node/.style={draw,circle}]
\def \n {5}
\def \radius {2cm}

  \node            (a) at (1.5,0)  { $v_{1}$ };
  \node            (b) at (3.0,0)  { $v_{2}$ };
  \node            (c) at (4.5,0)  { $v_{3}$ };
  \node            (d) at (6.0,0)  { $v_{4}$ };
  \node            (e) at (7.5,0)  { $v_{5}$ };
  \node            (f) at (9.0,0)  { $v_{6}$ };

  \draw[->] (f) edge [out=143,in=37] (a);
  \draw[->] (f) edge [out=145,in=35] (c);
  \draw[->] (d) edge [out=145,in=35] (a);
  \draw[->] (e) edge [out=210,in=330] (b);

\end{tikzpicture}
\label{fig:k6}
\end{center}

\begin{center} 
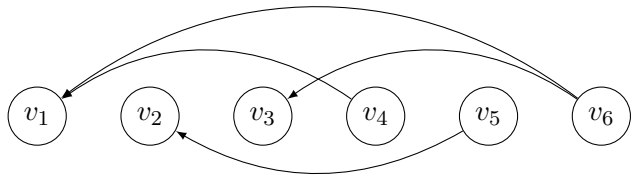
Fig.1 Tournament $K_{6}$. The only prime tournament on at most six vertices for which the conjecture is still open. Presented is the canonical ordering of its vertices. All edges that are not drawn are from left to right.
\end{center}

In this paper we prove the following:

\begin{theorem}
\label{six_vertex_theorem}
If $H$ is a six-vertex tournament not isomorphic to $K_{6}$ then it has the Erd\H{o}s-Hajnal property.
\end{theorem}

This reduces the six-vertex case to a single tournament.
The correctness of the conjecture for $K_{6}$ remains an open question.
Note that $K_{6}$ is a prime tournament. One can also check that $K_{6}$ does not have a galaxy ordering of vertices.
In fact the only ordering under which the graph of backward edges of $K_{6}$ is a forest is the canonical ordering presented in Fig.1.\\

We need to define two more special tournaments on six vertices that we denote by $L_{1}$ and $L_{2}$
and one special tournament on five vertices, denoted by $C_{5}$.

Tournament $C_{5}$ (see: Fig.2) is the unique tournament on five vertices such that each of its vertices has exactly two outneighbors 
and two inneighbors. Tournament $C_{5}$ is prime and one can check that it is not a galaxy.

\begin{center}
\begin{tikzpicture}[every path/.style={>=latex},every node/.style={draw,circle}]
\def \n {5}
\def \radius {2cm}

  \node[draw, circle] (a) at ({360/\n * (1 - 1)+18}:\radius) {$v_{1}$};
  \node[draw, circle] (b) at ({360/\n * (2 - 1)+18}:\radius) {$v_{2}$};
  \node[draw, circle] (c) at ({360/\n * (3 - 1)+18}:\radius) {$v_{3}$};
  \node[draw, circle] (d) at ({360/\n * (4 - 1)+18}:\radius) {$v_{4}$};
  \node[draw, circle] (e) at ({360/\n * (5 - 1)+18}:\radius) {$v_{5}$};

  \draw[->] (a) edge (b);
  \draw[->] (b) edge (c);
  \draw[->] (c) edge (d);
  \draw[->] (d) edge (e);
  \draw[->] (e) edge (a);
  \draw[->] (a) edge (c);
  \draw[->] (b) edge (d);
  \draw[->] (c) edge (e);
  \draw[->] (d) edge (a);
  \draw[->] (e) edge (b);

\end{tikzpicture}
\label{fig:c5}  
\end{center}
\begin{center} 
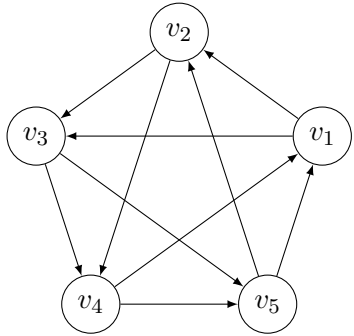
Fig.2 Tournament $C_{5}$ - the only prime five-vertex tournament that is not a galaxy.
\end{center}

Tournament $L_{1}$ is obtained from $C_{5}$ by adding one extra vertex and making it adjacent to exactly one vertex of $C_{5}$ (it does not matter to which one since all tournaments obtained by  procedure are isomorphic). Tournament $L_{2}$ is obtained from $C_{5}$ by adding one extra vertex and making it adjacent from $3$ vertices of $C_{5}$ that induce a cyclic triangle (again, it does not matter which cyclic triangle since all tournaments obtained by  this procedure are isomorphic).
Both tournaments are presented on Fig.3.

\begin{center}
\begin{tikzpicture}[every path/.style={>=latex},every node/.style={draw,circle}]
\def \n {5}
\def \radius {2cm}

  \node[draw, circle][shift={(0:-5cm)}] (a) at ({360/\n * (1 - 1)+18}:\radius) {$v_{1}$};
  \node[draw, circle][shift={(0:-5cm)}] (b) at ({360/\n * (2 - 1)+18}:\radius) {$v_{2}$};
  \node[draw, circle][shift={(0:-5cm)}] (c) at ({360/\n * (3 - 1)+18}:\radius) {$v_{3}$};
  \node[draw, circle][shift={(0:-5cm)}] (d) at ({360/\n * (4 - 1)+18}:\radius) {$v_{4}$};
  \node[draw, circle][shift={(0:-5cm)}] (e) at ({360/\n * (5 - 1)+18}:\radius) {$v_{5}$};
  \node[draw, circle][shift={(0:-10cm)}] (f) at ({360/\n * (5 - 1)+18}:\radius) {$v_{6}$};

  \node[draw, circle][shift={(0:-13cm)}] (aa) at ({360/\n * (1 - 1)+18}:\radius) {$v_{1}$};
  \node[draw, circle][shift={(0:-13cm)}] (bb) at ({360/\n * (2 - 1)+18}:\radius) {$v_{2}$};
  \node[draw, circle][shift={(0:-13cm)}] (cc) at ({360/\n * (3 - 1)+18}:\radius) {$v_{3}$};
  \node[draw, circle][shift={(0:-13cm)}] (dd) at ({360/\n * (4 - 1)+18}:\radius) {$v_{4}$};
  \node[draw, circle][shift={(0:-13cm)}] (ee) at ({360/\n * (5 - 1)+18}:\radius) {$v_{5}$};
  \node[draw, circle][shift={(0:-18cm)}] (ff) at ({360/\n * (5 - 1)+18}:\radius) {$v_{6}$};

  \draw[->] (a) edge (b);
  \draw[->] (b) edge (c);
  \draw[->] (c) edge (d);
  \draw[->] (d) edge (e);
  \draw[->] (e) edge (a);
  \draw[->] (a) edge (c);
  \draw[->] (b) edge (d);
  \draw[->] (c) edge (e);
  \draw[->] (d) edge (a);
  \draw[->] (e) edge (b);
  \draw[->] (f) edge  [out=340,in=200] (e);
  \draw[->] (d) edge  (f);
  \draw[->] (f) edge  (c);
  \draw[->] (a) edge  (f);
  \draw[->] (b) edge [out=185,in=80] (f);

  \draw[->] (aa) edge (bb);
  \draw[->] (bb) edge (cc);
  \draw[->] (cc) edge (dd);
  \draw[->] (dd) edge (ee);
  \draw[->] (ee) edge (aa);
  \draw[->] (aa) edge (cc);
  \draw[->] (bb) edge (dd);
  \draw[->] (cc) edge (ee);
  \draw[->] (dd) edge (aa);
  \draw[->] (ee) edge (bb);
  \draw[->] (ff) edge  [out=340,in=200] (ee);
  \draw[->] (dd) edge  (ff);
  \draw[->] (cc) edge  (ff);
  \draw[->] (aa) edge  (ff);
  \draw[->] (bb) edge [out=185,in=80] (ff);

\end{tikzpicture}
\end{center}

\begin{center} 
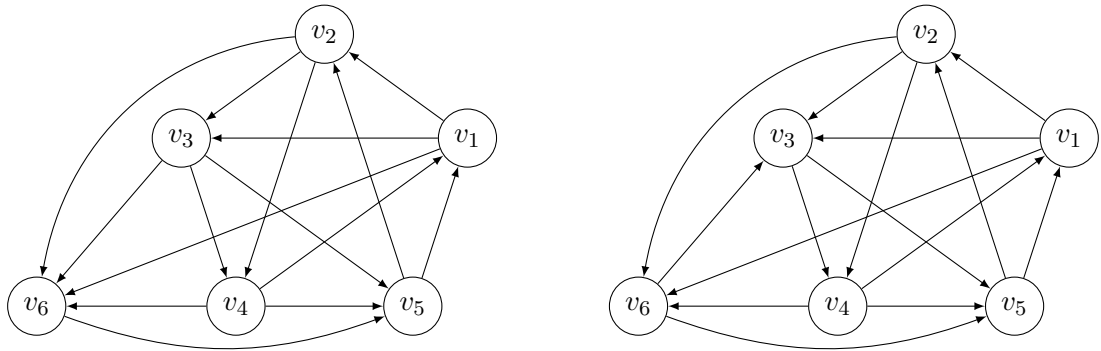
Fig.3 Tournament $L_{1}$ on the left and tournament $L_{2}$ on the right. Both are obtained from $C_{5}$ by adding one extra vertex.
\end{center}


This paper is organized as follows:
\begin{itemize}
\item in Section 2 we reduce the question about the correctness of the conjecture for six-vertex tournaments to three tournaments: $K_{6}, L_{1}, L_{2}$,
\item in Section 3 we introduce some tools to analyze tournaments $L_{1}$ and $L_{2}$,
\item in Section 4 we prove the conjecture for tournaments $L_{1}$ and $L_{2}$ and complete the proof of our main result.
\end{itemize}

\section{The landscape of six-vertex tournaments}

Our main result in this section is as follows:

\begin{theorem}
 \label{technical-theorem}
If $H$ is a six-vertex tournament not isomorphic to $K_{6}, L_{1}, L_{1}^{c}, L_{2}, L_{2}^{c}$ then $H$ satisfies the Erd\H{o}s-Hajnal conjecture.  
\end{theorem}

We will first prove a lemma describing the structure of
all six-vertex tournaments.

\begin{theorem}
 \label{technical-lemma}
Let $H$ be a six-vertex tournament. Then one of the following holds:
\begin{enumerate}
 \item $H$ is a galaxy, or
 \item there exists $v \in V(H)$, s.t. $H \setminus \{v\}$ is isomorphic to $C_{5}$ and $v$ has exactly one inneighbor or exactly one outneighbor in $H \setminus \{v\}$, or
 \item $H$ is not prime, or
 \item the vertices of $H$ or $H^{c}$ can be ordered as: $(a,b,c,d,e,f)$ such that the backward edges are: \\$(f,a), (e,a), (d,b), (f,c)$ (thus $H \setminus \{b\}$ or $H^{c} \setminus \{b\}$  is isomorphic to $C_{5}$
       and the outneighbors of $b$ form a cyclic triangle), or
 \item $H$ is isomorphic to $K_{6}$.
\end{enumerate}

\end{theorem}

\Proof
We may assume that $H$ is prime (for otherwise (3) holds), and so every vertex
of $H$ has at most four  inneighbors and at most four outneighbors.\\

\textbf{Case 1: some vertrex of $H$ has four outneighbors\\}
Suppose that $H$ has a vertex $v$ with $4$ outneighbors. Let $\{a,b,c,d\}$ be outneighbors of $v$ and denote by $u$ the remaining vertex.
Then $u$ is adjacent to $v$ and, since $H$ is prime, $u$  has at least one 
and at most $3$  outneighbors in $\{a,b,c,d\}$.

We call an ordering of the vertices of $H \setminus v$ {\em useful} if
it is   a galaxy ordering of $H \setminus v$, and no backward edge is 
incident with $u$. We observe that if $H \setminus v$ admits a useful ordering,
then adding $v$ at the start of this ordering produces a
galaxy ordering of $H$ (since $(u,v)$ is the only new backward edge, and no
other backward edge is incident with either $u$ or $v$), and  (1) holds. 
Thus we may assume that $H \setminus v$ admits no useful ordering.

Suppose first that $u$ has exactly three outneighbors in $\{a,b,c,d\}$, say $u$ 
is adjacent to  $a,b,c$ and from $d$. If $H|\{a,b,c\}$ is a transitive 
tournament (where $(a,b,c)$ is the transitive ordering, say), then $(d,u,a,b,c)$
is a useful ordering of $H \setminus v$, a contradiction. 
Therefore we may assume that $\{a,b,c\}$ induces a cyclic triangle. 
Without loss of generality we may assume that $(a,b), (b,c), (c,a)$ are edges.
Suppose first that $d$ has at most one inneighbor in $\{a,b,c\}$, say $b$ (without loss of generality) if one exists. But then $(d,u,a,b,c)$ is a useful
ordering of $H \setminus v$, a contradiction.
Thus $d$ has at least two inneighbors in $\{a,b,c\}$, i.e. $d$ has at most one outneighbor in $\{a,b,c\}$, say $b$ (without loss of generality) if one exists. But then $(u,v,a,b,c,d)$ is a galaxy ordering with backward edges: $(d,u), (c,a)$ and 
$(d,b)$ (if $b$ is an outneighbor of $d$), and so (1) holds.
We can thus assume that $u$ has at most two outneighbors in $\{a,b,c,d\}$. 

Next suppose that $u$ has exactly two outneighbors in $\{a,b,c,d\}$, say $u$ is adjacent from $a,b$ and to $c,d$.  Without loss of generality we  assume that $a$ is adjacent to $b$, and $c$ is adjacent to $d$. If there are at most $2$ edges from $\{c,d\}$ to $\{a,b\}$, then $(a,b,u,c,d)$ is a useful ordering of $H \setminus v$,
a contradiction. Thus we may assume that there are at least $3$ edges from $\{c,d\}$ to $\{a,b\}$. In other words, there is at most one edge from $\{a,b\}$ to $\{c,d\}$. If such an edge does not exist (i.e. $\{c,d\}$ is complete to $\{a,b\}$) then $(v,c,d,a,b,u)$ is a galaxy ordering of $H$, where each backward edge 
is incident with $u$, and (1) holds, so we may assume that  there is exactly 
one edge from $\{a,b\}$ to $\{c,d\}$.  We now check that in all cases
the theorem holds. If $a$ is adjacent to $d$ then $(v,c,a,d,b,u)$ is a galaxy 
ordering with all backward edges incident with $u$, and (1) holds.
If $b$ is adjacent to $c$ then $\{a,b,u,c,d\}$ induces a tournament isomorphic to $C_{5}$ and $v$ has a unique inneighbor in it, so (2) holds.
If $a$ is adjacent to $c$ then $\{u,v,d,a,c,b\}$ is a galaxy ordering with backward edges: $(a,u), (b,u), (c,d)$, and (1) holds. 
Finally, if $b$ is adjacent to $d$ then $(v,c,b,d,a,u)$ is a galaxy ordering with backward edges: $(u,v),(u,c),(u,d), (a,b)$, and again (1) holds.

Thus we may assume that $u$ has exactly one outneighbor in $\{a,b,c,d\}$, say
$a$. Let $(a^{'},b^{'},c^{'},d^{'})$
 be the ordering of $\{a,b,c,d\}$ in which $a$ has no backward edges, and where the number of backward edges is minimum subject to the previous constraint. Note that such an ordering is always a galaxy ordering. But then $(v,a^{'},b^{'},c^{'},d^{'},u)$ is also a galaxy ordering, and (1) holds. \\

We conclude that if some vertex in $H$ has $4$ outneighbors then the theorem holds.
Thus we can assume that every vertex of $H$ has at most three outneighbors. We can also conclude
that every vertex of $H$ has at most three inneighbors. The latter is true since the statement of the theorem
is invariant under reversing directions of all the edges of $H$. Indeed, after reversing all the edges the galaxy remains a galaxy, and  the property of being 
prime is also trivially invariant under this operation. Furthermore, both $C_{5}$ and $K_{6}$ are isomorphic to the tournaments obtained by reversing their edges. Therefore it remains to handle:\\

\textbf{Case 2: Every vertex has at most three outneighbors and at most three inneighbors\\}
Let us denote by $n_{3,2}$ the number of  vertices $v$ of $H$ such that $v$ has $3$ outneighbors and 
$2$ inneighbors. Similarly, let us denote by $n_{2,3}$ the number of vertices $v$ of $H$ such that $v$ has $3$ inneighbors and $2$ outneighbors. Then we have:
\begin{equation}
15=E(H) = 3n_{3,2} + 2n_{2,3} = 2n_{3,2} + 3n_{2,3}.
\end{equation}
Thus we have: $n_{3,2}=n_{2,3} = 3$.
Let $a,b,c$ be the vertices that have three outneighbors,  let  $x,y,z$ the remaining vertices.

Assume  first that $H|\{a,b,c\}$ is a transitive tournament, where $(a,b,c)$ 
(say) is a transitive ordering.
Then $c$ is complete to $\{x,y,z\}$ since, by definition, it has $3$ outneighbors, but it has no outneighbors in $\{a,b\}$. Similarly, vertex $b$ has exactly $2$ outneighbors in $\{x,y,z\}$ and without loss of generality we can assume that these are: $y$ and $z$. Vertex $a$ has exactly one outneighbor in $\{x,y,z\}$. 

Suppose first that $a$ is adjacent from $x$. Then, since $x$ has $2$ outneighbors and we already know that $x$ is adjacent to $a$ and $b$, we conclude that $x$ is adjacent from $y$ and $z$. Without loss of generality we can assume that $y$ is adjacent to $z$. If $a$ is adjacent to $y$ (and thus from $z$) then
$\{c,y\}$ is a homogeneous set and (3) holds. Thus we may assume that $a$ is adjacent to $z$ and from $y$.
But note that now $ H \setminus \{z\}$ is isomorphic to $C_{5}$ and $z$ has a unique outneighbor in $H \setminus \{z\}$, namely $x$. Thus (2) holds.
Therefore we may assume that $a$ is adjacent to $x$ and from $y$ and $z$.
Since $x$ has $2$ outneighbors, without loss of generality we can assume that $x$ is adjacent to $y$ and from $z$. Now, since $y$ has $2$ outneighbors, we can deduce that $y$ is adjacent to $z$
(this is true because the only outneighbor of $y$ in $\{a,b,c,x\}$ is $a$).
Now, $(a,c,x,b,y,z)$ is an ordering as in (4). This completes the case when
$H|{a,b,c}$ is a transitive tournament.

Thus we only need to consider the case when $\{a,b,c\}$ induces a cyclic 
triangle.
If $\{x,y,z\}$ induces a transitive tournament then we can reverse the edges of $H$ and repeat the analysis that we have just done for $\{a,b,c\}$. We can do it since, as we have already mentioned, the statement of the theorem is invariant under the operation of reversing all the edges of the tournament. Thus, without loss of generality, we can assume that both $\{x,y,z\}$ and $\{a,b,c\}$ induce cyclic triangles.
We may assume without loss of generality that $(x,y),(y,z),(z,x)$ and $(a,b),(b,c),(c,a)$ are edges.
Note that the edges from $\{x,y,z\}$ to $\{a,b,c\}$ form a matching. Indeed, each vertex of $\{x,y,z\}$ has exactly one outneighbor in $\{x,y,z\}$, therefore it has exactly one outneighbor in $\{a,b,c\}$ (since each vertex of $\{x,y,z\}$ has exactly $2$ outneighbors in $V(H)$), and each vertex from $\{a,b,c\}$ has exactly one inneighbor from $\{x,y,z\}$.

Without loss of generality we can assume that $x$ is adjacent to $a$. Assume first that $y$ is adjacent to $b$, and so $z$ is adjacent to $c$.
Now $(b,c,x,a,y,z)$ is a galaxy ordering with the backward edges: $(a,b)$, $(y,b)$, $(z,c)$, $(z,x)$.
Thus we may assume that $y$ is adjacent to $c$, and $z$ is adjacent to $b$.
But now $(b,c,x,a,y,z)$ is a canonical ordering of $K_{6}$, and (5) holds. That completes the proof of the lemma.
\bbox

We are now ready to prove Theorem \ref{technical-theorem}.

\Proof
We will use Lemma \ref{technical-lemma}. If outcome (1) holds then the result
follows from \ref{galaxy-theorem}. If outcome (2) holds then $H$ is isomorphic to one of the two tournaments: $L_{1}, L_{1}^{c}$. If outcomes (3) holds, then 
the result follows from \ref{nonprime} and  \ref{five-vertex-theorem}.
Finally, if outcome (4) holds then $H$ is isomorphic to $L_{2}$ or $L_{2}^{c}$.
This completes the proof of Theorem \ref{technical-theorem}.
\bbox

\section{Regularity tools}

In this section we will introduce some regularity tools that will be very useful later on to prove the conjecture for $L_{1}$ and $L_{2}$.

Denote by $tr(T)$ the largest size of the transitive subtournament of $T$. 
For $X \subseteq V(T)$, write $tr(X)$ for $tr(T|X)$.
Let $X,Y \subseteq V(T)$ be disjoint. Denote by $e_{X,Y}$ the number of directed edges $(x,y)$, where $x \in X$ and $y \in Y$.
The \textit{directed density from X to Y} is defined as 
$d(X,Y)=\frac{e_{X,Y}}{|X||Y|}.$ 

We call a tournament $T$  \textit{$\epsilon$-critical} for $\epsilon>0$ if $tr(T) < |T|^{\epsilon}$ but for every proper subtournament $S$ of $T$ we have: $tr(S) \geq |S|^{\epsilon}$. Next we list some properties of $\epsilon$-critical tournaments that we borrow from \cite{chorochudber}.

\begin{theorem}
\label{remarktheorem}
For every $N>0$ there exists $\epsilon(N)>0$ such that for every $0<\epsilon<\epsilon(N)$ every $\epsilon$-critical tournament $T$ satisfies $|T| \geq N$.
\end{theorem}

\Proof
Since every tournament contains a transitive subtournament of order $2$ so it suffices to take $\epsilon(N)=\log_{N}(2)$. \bbox 

\begin{theorem}
\label{firsttechnicallemma}
Let $T$ be an $\epsilon$-critical tournament with $|T|=n$ and $\epsilon,c,f>0$ be constants such that $\epsilon <\log_{c}(1-f)$.
Then for every $A \subseteq V(T)$ with $|A| \geq cn$ and every transitive subtournament $G$ of $T$ with $|G| \geq f \cdot tr(T)$ we have: $A$ is not complete from $V(G)$ and $A$ is not complete to $V(G)$.
\end{theorem}

\Proof
Assume otherwise. Let $A_T$ be a transitive subtournament in $T|A$ of size 
$tr(A)$.
Then $|A_{T}| \geq (cn)^{\epsilon}$. Now we can merge $A_{T}$ with $G$ to obtain a transitive subtournament of size at least $(cn)^{\epsilon}+f tr(T)$. From the definition of $tr(T)$ we have $(cn)^{\epsilon}+f tr(T) \leq tr(T)$. So $c^{\epsilon}n^{\epsilon} \leq (1-f) tr(T)$, and in particular $c^{\epsilon}n^{\epsilon} < (1-f)n^{\epsilon}$. But this contradicts the fact that $\epsilon <\log_{c}(1-f)$. \bbox

\begin{theorem}
\label{veryeasylemma}
Let  $T$ be an $\epsilon$-critical tournament with $|T|=n$ and $\epsilon,c>0$ be constants such that $\epsilon<\log_{\frac{c}{2}}(\frac{1}{2})$.
Then for every two disjoint subsets $X,Y \subseteq V(T)$ with $|X| \geq cn$, $|Y| \geq cn$ 
there exist an integer 
$k \geq \frac{cn}{2}$ 
and vertices $x_{1},...,x_{k} \in X$ and $y_{1},...,y_{k} \in Y$ such that  $y_{i}$ is adjacent to 
$x_{i}$ for $i=1,...,k$. 
\end{theorem}

\Proof
Assume otherwise. Write $m= \lfloor \frac{cn}{2} \rfloor$.
Consider the bipartite graph $G$ with bipartition $(X,Y)$ where $\{x,y\} \in E(G)$ if $(y,x) \in V(T)$.
Then we know that $G$ has no matching of size $m$.
By K\"{o}nig's Theorem (see \cite{diestel}) there exists $C \subseteq V(G)$
with $|C|<m$, such that every edge of $G$ has an end in $C$. Write
$C \cap X=C_{X}$ and $C \cap Y = C_{Y}$. We have $|C_{X}| \leq \frac{|X|}{2}$ and $|C_{Y}| \leq \frac{|Y|}{2}$. Therefore $|X \backslash C_{X}| \geq \frac{|X|}{2}$ and $|Y \backslash C_{Y}| \geq \frac{|Y|}{2}$, and
by the definition of $C$ and $G$, we know that $X \backslash C_{X}$ is complete to $Y \backslash C_{Y}$. Denote by $T_{1}$ a transitive subtournament of size
$tr(T|(X \backslash C_{X}))$ in $T|(X \backslash C_{X})$. Denote by $T_{2}$ a transitive subtournament of size $tr(T|(Y \backslash C_{Y}))$ in 
$T|(Y \backslash C_{Y})$.
From the $\epsilon$-criticality of $T$ and since $|X \backslash C_{X}| \geq \frac{cn}{2}$, $|Y \backslash C_{Y}| \geq \frac{cn}{2}$, we also have: $|T_{1}| \geq (\frac{cn}{2})^{\epsilon}$, $|T_{2}| \geq (\frac{cn}{2})^{\epsilon}$.
We can merge $T_{1}$ and $T_{2}$ to obtain bigger transitive tournament $T_{3}$ with $|T_{3}| \geq 2(\frac{c}{2})^{\epsilon} n^{\epsilon}$.
Therefore, since $T$ is $\epsilon$-critical, we have: $2(\frac{c}{2})^{\epsilon} < 1$. But this contradicts the condition $\epsilon<\log_{\frac{c}{2}}(\frac{1}{2})$.
          
\bbox

Next we introduce one more structure that will be crucial to prove the conjecture for $L_{1}$
and $L_{2}$. Again, its definition can be found in \cite{chorochudber}, but we give it again for the
reader's convenience.

Let $c>0$, $0<\lambda<1$ be constants, and let $w$ be a $\{0,1\}$-vector of length $|w|$. Let $T$ be a tournament with $|T|=n$. A sequence of disjoint subsets 
$(S_{1},S_{2},...,S_{|w|})$ of $V(T)$ is a $(c, \lambda, w)$-{\em structure} if
\begin{itemize} 
\item whenever $w_i=0$ we have $|S_{i}| \geq cn$ (we say that $S_i$ is a {\em linear set}) 
\item whenever $w_i=1$ the set $T|S_{i}$ is transitive and $|S_{i}| \geq c \cdot tr(T)$ (we say that $T_i$ is a {\em transitive set})
\item $d(S_{i},S_{j}) \geq 1 - \lambda$ for all $1 \leq i < j \leq |w|$.
\end{itemize} 

The following was proved in \cite{chorochudber}:

\begin{theorem}
\label{thirdtechnicallemma}
Let $S$ be a tournament, let $w$ be a $\{0,1\}$-vector, and let 
$0 < \lambda_{0} < \frac{1}{2}$ be a constant. Then there exist $\epsilon_{0},c_{0}>0$ such that for every $0<\epsilon<\epsilon_{0}$, every $S$-free $\epsilon$-critical tournament contains a $(c_{0}, \lambda_{0}, w)$-structure.
\end{theorem}

\begin{center}
\begin{tikzpicture}[every path/.style={>=latex},every node/.style={draw,circle}]
\def \n {5}
\def \radius {3cm}

  \node[input]            (a) at (3.3,0)  { $A_{1}$ };
  \node[matrx]            (b) at (6.6,0)  { $T_{1}$ };
  \node[input]            (c) at (9.9,0)  { $A_{2}$ };
  \node[input]            (d) at (13.5,0)  { $A_{3}$ };
  \node[matrx]            (e) at (16.8,0)  { $T_{2}$ };

  \draw[vecArrow] (a) to (b);
  \draw[innerWhite] (a) to (b);

  \draw[vecArrow] (b) to (c);
  \draw[innerWhite] (b) to (c);

  \draw[vecArrow] (c) to (d);
  \draw[innerWhite] (c) to (d);

  \draw[vecArrow] (d) to (e);
  \draw[innerWhite] (d) to (e);

\end{tikzpicture}
\end{center}

\begin{center} 
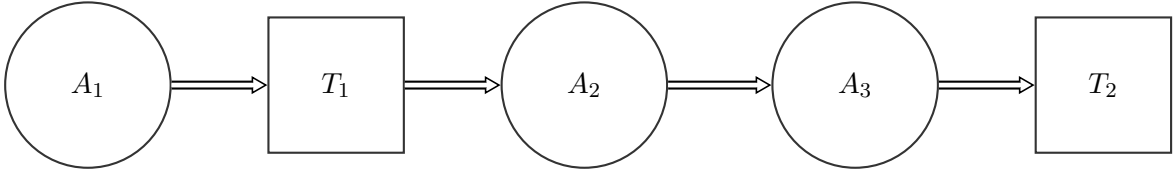
Fig.4 Schematic representation of the $(c,\lambda,w)$-structure. This structure consists of three linear sets: $A_{1},A_{2},A_{3}$ and two transitive sets: $T_{1}$ and $T_{2}$. The arrows indicate the orientation of most
of the edges going between different elements of the $(c,\lambda,w)$-structure. Each $T_{i}$ satisfies: $|T_{i}| \geq c \cdot tr(T)$ and each $A_{i}$ satisfies: $|A_{i}| \geq c \cdot n$, where $n=|T|$. We have here: $w=(0,1,0,0,1)$. 
\end{center}

We say that a $(c,\lambda,w)$-structure is \textit{smooth} if the last condition of the definition of the
$(c,\lambda,w)$-structure is satisfied in a stronger form, namely we have: $d(\{v\}, S_{j}) \geq 1 - \lambda$ for $v \in S_{i}$ and $d(S_{i},\{v\}) \geq 1 - \lambda$ for $v \in S_{j}$, $i<j$.

Theorem \ref{thirdtechnicallemma} leads to the following  conclusion: 

\begin{theorem}
\label{smooththeorem}
Let $S$ be a tournament, let $w$ be a $\{0,1\}$-vector, and let 
$0 < \lambda_{1} < \frac{1}{2}$ be a constant. Then there exist $\epsilon_{1},c_{1}>0$ such that for every $0<\epsilon<\epsilon_{1}$, every $S$-free $\epsilon$-critical tournament contains a smooth $(c_{1}, \lambda_{1}, w)$-structure.
\end{theorem}

\Proof

By Theorem \ref{thirdtechnicallemma}, there exist 
$\epsilon_{0},c_{0}>0$ such that for every $0<\epsilon<\epsilon_{0}$, every $S$-free $\epsilon$-critical tournament contains a $(c_{0}, \lambda_{0}, w)$-structure.
Denote this structure by $(A_{1},...,A_{k})$. Let $M$ be a positive constant. 
For an ordered pair $(i,j)$, where $i,j \in \{1,...,k\}$ and $i \neq j$ let $Bad^{M}(i,j)$ be the set of these vertices $v \in A_{i}$ such that 
\begin{itemize}
 \item $v$ is adjacent from more than $M \lambda_{0} |A_{j}|$ vertices of $A_{j}$ if $i < j$ and
 \item $v$ is  adjacent to more than $M \lambda_{0} |A_{j}|$ vertices of $A_{j}$ if $i > j$.
\end{itemize}
Note first that $|Bad^{M}(i,j)| \leq \frac{|A_{i}|}{M}$. Indeed, otherwise by the definition of $Bad^{M}(i,j)$, the number of backward edges between $A_{i}$ and $A_{j}$
is more than $\lambda_{0}|A_{i}||A_{j}|$ which contradicts the fact that $d(A_{\min(i,j)},A_{\max(i,j)}) \geq 1 - \lambda_{0}$.
Now let  $A_{i}^{M} = A_{i} \setminus \bigcup_{j \in \{1,...,k\}, j \neq i} Bad^{M}(i,j)$.
From the fact that $|Bad^{M}(i,j)| \leq \frac{|A_{i}|}{M}$, we get $|A_{i}^{M}| \geq (1-\frac{k-1}{M})|A_{i}|$. Now take $M=2k$.
Then we obtain $|A_{i}^{M}| \geq \frac{|A_{i}|}{2}$.
Consider the sequence $(A_{1}^{M},...,A_{k}^{M})$. Take a pair $\{i,j\}$, where $i,j \in \{1,...,k\}$ and $i < j$.
Note that by the definition of $A_{i}^{M}$, we know that every vertex $v \in A_{i}^{M}$ is adjacent from at most $M\lambda_{0}|A_{j}|$ vertices of
$A_{j}^{M}$. For $M=2k$, since $|A_{j}^{M}| \geq \frac{|A_{j}|}{2}$, we obtain: every vertex $v \in A_{i}^{M}$ is adjacent from at most
$2M\lambda_{0}|A_{j}^{M}|$ vertices of $A_{j}^{M}$. Similarly, we get: every vertex $v \in A_{j}^{M}$ is adjacent to at most
$2M\lambda_{0}|A_{i}^{M}|$ vertices of $A_{i}^{M}$. Consequently,  $(A_{1}^{M},...,A_{k}^{M})$ is a smooth $(\frac{c_{0}}{2},2M\lambda_{0},w)$-structure.
Thus taking: $\lambda_{0} = \frac{\lambda_{1}}{4k}$ and $c_{1} = \frac{c_{0}}{2}$, we complete the proof.
\bbox

\section{The Erd\H{o}s-Hajnal conjecture holds for $L_{1}$ and $L_{2}$}

We are ready to prove that both $L_{1}$ and $L_{2}$ satisfy the conjecture.
We will use two special orderings of the vertices of $L_{1}$ and two special 
orderings of the vertices of $L_{2}$.

\begin{center}
\begin{tikzpicture}[every path/.style={>=latex},every node/.style={draw,circle}]
\def \n {5}
\def \radius {3cm}

  \node            (a) at (1.2,0)  { $v_{3}$ };
  \node            (b) at (2.4,0)  { $v_{4}$ };
  \node            (c) at (3.6,0)  { $v_{5}$ };
  \node            (d) at (4.8,0)  { $v_{1}$ };
  \node            (e) at (6.0,0)  { $v_{2}$ };
  \node            (f) at (7.2,0)  { $v_{6}$ };

  \node            (g) at (10.2-1,0)  { $v_{2}$ };
  \node            (h) at (11.4-1,0)  { $v_{4}$ };
  \node            (i) at (12.6-1,0)  { $v_{1}$ };
  \node            (j) at (13.8-1,0)  { $v_{3}$};
  \node            (k) at (15.0-1,0)  { $v_{6}$ };
  \node            (l) at (16.2-1,0)  { $v_{5}$ };

  \draw[->] (e) edge [out=150,in=30] (a);
  \draw[->] (d) edge [out=210,in=330] (a);
  \draw[->] (e) edge [out=210,in=330] (b);
  \draw[->] (f) edge [out=150,in=30] (c);

  \draw[->] (l) edge [out=150,in=30] (i);
  \draw[->] (i) edge [out=150,in=30] (g);
  \draw[->] (l) edge [out=150,in=30] (g);
  \draw[->] (j) edge [out=215,in=325] (h);



\end{tikzpicture}
\end{center}

\begin{center} 
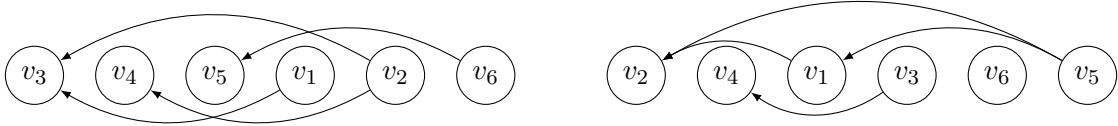
Fig.5 Two crucial orderings of the vertices of $L_{1}$. The left one is the forest ordering and the right one is the cyclic ordering. Notice that neither of them is a galaxy ordering.
\end{center}

The first ordering of the vertices of $L_{1}$ is as follows: $(v_{3},v_{4},v_{5},v_{1},v_{2},v_{6})$, where the set of backward edges is: $\{(v_{1},v_{3}),(v_{2},v_{4}),(v_{2},v_{3}),(v_{6},v_{5})\}$. We call it the \textit{forest ordering of $L_{1}$} since under this ordering the graph of backward edges is a forest.
The second ordering of the vertices of $L_{1}$ is as follows: $(v_{2},v_{4},v_{1},v_{3},v_{6},v_{5})$, where the set of backward edges is: $\{(v_{1},v_{2}),(v_{5},v_{1}),(v_{5},v_{2}),(v_{3},v_{4})\}$.
We call it the \textit{cyclic ordering of $L_{1}$}.

\begin{center}
\begin{tikzpicture}[every path/.style={>=latex},every node/.style={draw,circle}]
\def \n {5}
\def \radius {3cm}

  \node            (a) at (1.2,0)  { $v_{1}$ };
  \node            (b) at (2.4,0)  { $v_{2}$ };
  \node            (c) at (3.6,0)  { $v_{3}$ };
  \node            (d) at (4.8,0)  { $v_{4}$ };
  \node            (e) at (6.0,0)  { $v_{6}$ };
  \node            (f) at (7.2,0)  { $v_{5}$ };

  \node            (g) at (10.2-1,0)  { $v_{2}$ };
  \node            (h) at (11.4-1,0)  { $v_{4}$ };
  \node            (i) at (12.6-1,0)  { $v_{1}$ };
  \node            (j) at (13.8-1,0)  { $v_{6}$};
  \node            (k) at (15.0-1,0)  { $v_{3}$ };
  \node            (l) at (16.2-1,0)  { $v_{5}$ };

  \draw[->] (f) edge [out=150,in=30] (a);
  \draw[->] (d) edge [out=210,in=330] (a);
  \draw[->] (f) edge [out=210,in=330] (b);
  \draw[->] (e) edge [out=150,in=30] (c);

  \draw[->] (l) edge [out=150,in=30] (i);
  \draw[->] (i) edge [out=150,in=30] (g);
  \draw[->] (l) edge [out=150,in=30] (g);
  \draw[->] (k) edge [out=215,in=325] (h);



\end{tikzpicture}
\end{center}

\begin{center} 
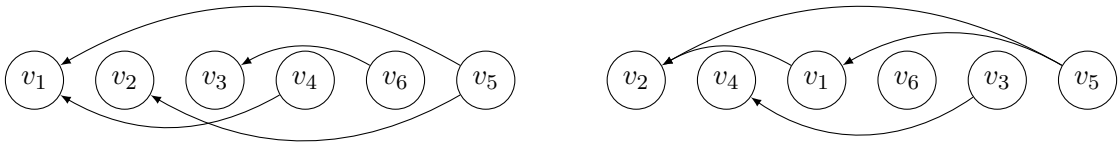
Fig.6 Two crucial orderings of vertices of $L_{2}$. The left one is the forest ordering and the right one is the cyclic ordering. Notice that neither of them is a galaxy ordering.
\end{center}

The first ordering of the vertices of $L_{2}$ is as follows: $(v_{1},v_{2},v_{3},v_{4},v_{6},v_{5})$, where the set of backward edges is: $\{(v_{4},v_{1}),(v_{5},v_{2}),(v_{5},v_{1}),(v_{6},v_{3})\}$. We call it the \textit{forest ordering of $L_{2}$}.
The second ordering of the vertices of $L_{2}$ is as follows: $(v_{2},v_{4},v_{1},v_{6},v_{3},v_{5})$, where the set of backward edges is: $\{(v_{1},v_{2}),(v_{5},v_{1}),(v_{5},v_{2}),(v_{3},v_{4})\}$.
We call it the \textit{cyclic ordering of $L_{2}$}.

\begin{theorem}
\label{l2-theorem}
Tournament $L_{2}$ satisfies the Erd\H{o}s-Hajnal conjecture.
\end{theorem}

\Proof
We will prove that every $L_{2}$-free tournament $T$ on $n$ vertices contains a transitive subtournament
of size at least $n^{\epsilon}$ for $\epsilon > 0$ small enough.
Assume for a contradiction that this is not the case and let
$T$ be the smallest $L_{2}$-free $\epsilon$-critical tournament. By Theorem \ref{remarktheorem} we may assume that $|T|$ is large enough. We will get a contradiction, proving that $T$ contains a transitive subtournament of order $n^{\epsilon}$. 
By Theorem \ref{smooththeorem} we  extract from $T$ a smooth $(c_{0}(\lambda_{0}), \lambda_{0}, w)$-structure $\chi_{0}=(A_{1},A_{2},T_{0},A_{3},A_{4},A_{5})$,
where $w=(0,0,1,0,0,0)$ and $\lambda_{0} > 0$ is an arbitrary positive number.
We will fix $\lambda_{0}$ to be small enough. 
We then take an arbitrary subset $S$ of $T_{0}$ such that $|S|$ is divisible by $3$ and $|S|$ is of maximum size. Notice that $|S| \geq |T_{0}|-2$.
Since  $|T_{0}| \geq c_{0}(\lambda_{0}) tr(T)$ and $|T|$ is large, it follows
that $|T_{0}| \geq 4$, and so $|S| \geq \frac{|T_{0}|}{2}$.
Now take the sequence $\chi=(A_{1},A_{2},S,A_{3},A_{4},A_{5})$.
Since $(A_{1},A_{2},T_{0},A_{3},A_{4},A_{5})$ is a a smooth $(c_{0}(\lambda_{0}), \lambda_{0}, w)$-structure
and  $S$ is a subset of $T_{0}$ of size $|S| \geq \frac{|T_{0}|}{2}$, we get that $(A_{1},A_{2},S,A_{3},A_{4},A_{5})$
is a smooth $(c(\lambda), \lambda, w)$-structure for $\lambda = 2\lambda_{0}$ and $c(\lambda) = \frac{c_{0}(\lambda_{0})}{2} = \frac{c_{0}(\frac{\lambda}{2})}{2}$.
We partition $S$ into three subsets: the set of first $\frac{|S|}{3}$ vertices called $T_{1}$, the set of next $\frac{|S|}{3}$ vertices called $T_{2}$ and the remaining part called $T_{3}$ (here we refer to the transitive ordering of $S$).
By Theorem \ref{veryeasylemma} we may assume that there exist $x_{1},...,x_{k} \in A_{1}$
and $y_{1},...,y_{k} \in A_{5}$ such that $k \geq \frac{cn}{2}$ and $(y_{i},x_{i})$ is an edge for $i=1,...,k$. Denote $X=\{x_{1},...,x_{k}\}$, $Y=\{y_{1},...,y_{k}\}$.
Let $X_{wrong}$ be the set of vertices of $X$ that are complete to 
$T_{3}$, and let $Y_{wrong}$ the set of vertices of of $Y$ that are complete 
from  $T_{1}$. Assume first that $|X_{wrong}| \geq \frac{k}{3}$.
But $X_{wrong}$ is complete to $T_{3}$, $|X_{wrong}| \geq \frac{c}{6}n$,  and $|T_{3}| \geq \frac{c}{3}tr(T)$, which contradicts Theorem \ref{firsttechnicallemma} if $\epsilon < \log_{\frac{c}{6}}(1-\frac{c}{3})$. We get a
similar contradiction if $|Y_{wrong}| \geq \frac{k}{3}$.
Therefore $|X_{wrong}| < \frac{k}{3}$ and $|Y_{wrong}| < \frac{k}{3}$. Write $\mathcal{I} = \{i \in \{1,...,k\}\:x_{i} \notin X_{wrong} \land y_{i} \notin Y_{wrong}\}$. We have: $|\mathcal{I}| > \frac{k}{3}$, and in particular $\mathcal{I} \neq \emptyset$. Fix $j \in \mathcal{I}$. Let $u \in T_1$ be an outneighbor of $y_{j}$,
and let $v \in T_3$ be  an inneighbor of $x_{j}$.  Note that since $u \in T_{1}$ and $v \in T_{3}$, $(u,v)$ is an edge.

Assume first that both $(x_{j},u)$ and $(v,y_{j})$ are edges.
Let $T_{2}^{*}$ be the set of vertices of $T_2$ that are  outneighbors of 
$x_{j}$ and inneighbors of $y_{j}$. 
From the fact that $\chi$ is smooth, we get: $|T_{2}^{*}| \geq |T_{2}|-2\lambda|T| \geq \frac{c}{3}(1-6\lambda)tr(T) \geq \frac{c}{6}tr(T)$ if we take $\lambda \leq \frac{1}{12}$.
Let $A_{3}^{*}$ be the set of vertices of $A_{3}$ that are outneighbors of $x_{j}$, $u$ and $v$, and inneighbors of $y_{j}$. Again, from the fact that $\chi$ is smooth, we get:
$|A_{3}|^{*} \geq |A_{3}|(1-4\lambda) \geq \frac{c}{2}n$ for $\lambda \leq \frac{1}{8}$. 
Now, if $\epsilon < \log_{\frac{c}{2}}(1-\frac{c}{6})$,  by Theorem 
\ref{smooththeorem} there exists  $z \in A_{3}^{*}$  and $w \in T_{2}^{*}$ such 
that $(z,w)$ is an edge, and so  $(x_{j},u,w,v,z,y_{j})$ is the forest ordering 
of  $L_{2}$, a contradiction.

Thus either $(u,x_{j})$ is an edge or $(y_{j},v)$ is an edge.
Assume that the former holds (if the latter holds, the argument is similar, and
we omit it). Let $A_{2}^{*}$ be the set of vertices of  $A_{2}$ that are outneighbors of $x_{j}$ and inneighbors of $u$ and $y_{j}$. From the fact that $\chi$ is smooth, we get: $|A_{2}^{*}| \geq |A_{2}|(1-3\lambda) \geq \frac{c}{2}n$ for $\lambda \leq \frac{1}{6}$. Let $A_{4}^{*}$  be the set of vertices of $A_{4}$ that are outneighbors of $x_{j}$ and $u$, and inneighbors of $y_{j}$. From the fact that $\chi$ is smooth, we get: $|A_{4}^{*}| \geq |A_{4}|(1-3\lambda) \geq \frac{c}{2}n$ for $\lambda \leq \frac{1}{6}$.
Now, if  $\epsilon < \log_{\frac{c}{4}}(\frac{1}{2})$, Theorem \ref{veryeasylemma} implies that there exist
$z \in A_{4}^{*}$  and $ w \in A_{2}^{*}$ such that $(z,w)$ is an edge.
Let $A_{3}^{*}$ be the set of vertices of $A_{3}$ that are outneighbors of $x_{j},w,u$, and inneighbors of $z,y_{j}$. From the fact that $\chi$ is smooth, we get:
$|A_{3}^{*}| \geq |A_{3}|(1-5\lambda) \geq \frac{c}{2}n$ for $\lambda < \frac{1}{10}$.
In particular, $A_{3}^{*}$ is nonempty. Let $s \in A_{3}^{*}$.
Now $(x_{j},w,u,s,z,y_{j})$ is the cyclic ordering of $L_{2}$, again a contradiction. This completes the proof.
\bbox

\begin{theorem}
\label{l1-theorem}
Tournament $L_{1}$ satisfies the Erd\H{o}s-Hajnal conjecture.
\end{theorem}

\Proof
The proof goes along the same line as the proof of the previous theorem.

Again we take an $\epsilon$-critical tournament $T$ that this time is $L_{1}$-free, and 
get a contradiction for $\epsilon>0$ small enough.
By Theorem \ref{smooththeorem} we  extract from $T$ a smooth $(c_{0}(\lambda_{0}), \lambda_{0}, w)$-structure $\chi_{0}=(A_{1},A_{2},T_{0},A_{3},A_{4},A_{5}, A_{6})$,
where $w=(0,0,1,0,0,0,0)$ and $\lambda_{0} > 0$ is an arbitrary positive number.
We will fix $\lambda_{0}$ to be small enough.
As in the previous proof, we use $\chi_{0}$ to construct a $(c(\lambda),\lambda,w)$-structure $\chi=(A_{1},A_{2},S,A_{3},A_{4},A_{5}, A_{6})$, where
$|S|$ is divisible by $3$.
We partition $S$ into three subsets: the set of first $\frac{|S|}{3}$ vertices called $T_{1}$, the set of next $\frac{|S|}{3}$ vertices called $T_{2}$ and the remaining part 
called $T_{3}$.

As in the previous proof, we may assume that there exist $x_{j} \in A_{1},y_{j} \in  A_{5}$ 
such that $(y_{j},x_{j})$ is an edge, $y_{j}$ has an outneighbor $u$ in $T_{1}$, and $x_{j}$ has an inneighbor $v$ in  $T_{3}$.

Assume first that both $(x_{j},u)$ and $(v,y_{j})$ are edges.
Now denote by $T_{2}^{*}$ the set of vertices of $T_{2}$ that are outneighbors of $x_{j}$, and inneighbors of $y_{j}$. From the fact that $\chi$ is smooth, we get: $|T_{2}^{*}| \geq |T_{2}|-2\lambda|T| \geq \frac{c}{3}(1-6\lambda)tr(T) \geq \frac{c}{6}tr(T)$ if we take $\lambda \leq \frac{1}{12}$.
Let us also denote by $A_{6}^{*}$ the set of vertices of  $A_{6}$ that are outneighbors of $x_{j}$, $u$, $v$ and $y_{j}$. Again, from the fact that $\chi$ is smooth, we get:
$|A_{6}|^{*} \geq |A_{6}|(1-4\lambda) \geq \frac{c}{2}n$ for $\lambda \leq \frac{1}{8}$. 
Now, if $\epsilon < \log_{\frac{c}{2}}(1-\frac{c}{6})$, 
Theorem \ref{smooththeorem} implies that there  exist $z \in A_{6}^{*}$ 
and  $w \in T_{2}^{*}$ such that $(z,w)$ is an edge, and so 
$(x_{j},u,w,v,y_{j},z)$  is the forest ordering of $L_{1}$, a contradiction.

Thus either $(u,x_{j})$ is an edge or $(y_{j},v)$ is an edge.
We assume that the former holds (if the latter holds, the argument is similar
and we omit it).
Let  $A_{2}^{*}$ be the set of vertices of  $A_{2}$ that are outneighbors of $x_{j}$ and inneighbors of $u$ and $y_{j}$. From the fact that $\chi$ is smooth, we get: $|A_{2}^{*}| \geq |A_{2}|(1-3\lambda) \geq \frac{c}{2}n$ for $\lambda \leq \frac{1}{6}$. Let $A_{3}^{*}$ be the set of vertices of  $A_{3}$ that are outneighbors of $x_{j}$ and $u$, and inneighbors of $y_{j}$. From the fact that $\chi$ is smooth, we get: $|A_{3}^{*}| \geq |A_{3}|(1-3\lambda) \geq \frac{c}{2}n$ for $\lambda \leq \frac{1}{6}$.
Now, if   $\epsilon < \log_{\frac{c}{4}}(\frac{1}{2})$, 
Theorem \ref{veryeasylemma} implies that there exist 
$z \in A_{3}^{*}$  and $w \in A_{2}^{*}$ such that $(z,w)$ is an edge.
Denote by $A_{4}^{*}$ the set of vertices of $A_{4}$ that are outneighbors of $x_{j},w,u,z$, and inneighbors of $y_{j}$. From the fact that $\chi$ is smooth, we get:
$|A_{4}^{*}| \geq |A_{4}|(1-5\lambda) \geq \frac{c}{2}n$ for $\lambda < \frac{1}{10}$.
In particular, $A_{4}^{*}$ is nonempty. Let $s \in A_{4}^{*}$.
Now $(x_{j},w,u,z,s,y_{j})$ is a cyclic ordering of $L_{1}$, again a 
contradiction.
This completes the proof.
\bbox

We are now ready to finish the proof of Theorem \ref{six_vertex_theorem}.

\Proof
By Theorem \ref{technical-theorem}, it suffices to prove the conjecture for $L_{1},L_{1}^{c},L_{2}$ and $L_{2}^{c}$. We have just proved the conjecture for $L_{1}$ and $L_{2}$. Thus obviously $L_{1}^{c}$ and $L_{2}^{c}$ also satisfy the conjecture. This completes the proof.
\bbox


\end{document}